\documentclass[11pt]{article}
\usepackage{graphicx} 
\usepackage{latexsym} 
\usepackage{amsmath}
\usepackage{amssymb} 
\usepackage{verbatim} 
\usepackage{alltt} 
\begin{document}

\title{ 
A case of combination of evidence in the Dempster-Shafer theory 
inconsistent with evaluation of probabilities}

\author{Andrzej K. Brodzik and Robert H. Enders \\ \\ \\
The MITRE Corporation \\
202 Burlington Road \\
Bedford, MA 01730, USA \\
email: abrodzik@mitre.org \\
} 
\date{}
 
\maketitle
 
\baselineskip 4.30mm

\begin{abstract}
\noindent
The Dempster-Shafer theory of evidence accumulation is one of the main tools for combining data obtained from multiple sources.
In this paper a special case of combination of two bodies of evidence 
with non-zero conflict coefficient is considered.
It is shown that application of the Dempster-Shafer rule of combination in this case 
leads to an evaluation of masses of the combined bodies 
that is different from the evaluation of the corresponding probabilities obtained by application of the 
law of total probability. This finding supports the view that
probabilistic interpretation of results of the Dempster-Shafer analysis in the general case is not appropriate. 
\end{abstract}

\noindent

\section{Introduction}

\vspace{2mm}

\noindent
{\it The greatest enemy of knowledge is not ignorance, it is the illusion of knowledge} 
\begin{flushright}						
\ \ \ \ \ \ \ \ \hspace{40mm} \ \ \ \ \ \ \ \ \ Stephen Hawking \\ 
\end{flushright}

\vspace{1mm}

\noindent
Data fusion is the first stage of a complex decision making process. 
This stage usually involves combination of uncertain, incomplete, and/or difficult-to-compare information obtained from multiple sources. 
The principal goals of this task are to decrease uncertainty associated with individual measurements 
and to permit identification of the most likely alternative. 
One of the main tools for data fusion is the Dempster-Shafer (DS) theory of evidence accumulation \cite{dempster}, \cite{ds}. 
The DS theory has been used in many areas of science and engineering, including sensor fusion, medical diagnostics, 
image processing, biometrics, and decision support. A review of some of these applications is given in \cite{app}.

In this work we focus on certain key aspects of the relationship between the DS and probability theories, following up on prior work on this subject. 
Recently, an investigation of the algebraic structure of the DS set of mass assignments was undertaken and it was shown that this set
can be mapped onto the set of probabilities by a semigroup homomorphism \cite{brodzik}.
In particular, it was demonstrated that the combination of mass from the singleton DS set is in the abstract-algebraic sense equivalent to 
the combination of probabilities.
While a naive interpretation of this result supports the proposition that the DS theory in the non-singleton case is, 
in a certain sense, a generalization of probability theory \cite{bp}, there is an evidence that 
the two data fusion approaches differ in several respects that seem to defy such interpretation 
\cite{dia}.

Recently, several 
cases where the results of DS evidence accumulation might present interpretation difficulties from the probability theory standpoint have been observed. 
In one of these cases two bodies of evidence with mass assignments 0.99, 0.00, 0.01 and 0.00, 0.99, 0.01 are combined.
This results in the masses associated with the decision set 0.00, 0.00, 1.00 - an outcome that is deemed undesirable \cite{zadeh}.
This result occurs due to strongly contradictory beliefs about the first two elements. 
The problem can be relieved to some extent by replacing the zero mass assignments with appropriately small but non-zero values,
but it is not clear that an arbitrary resolution of such contradictions is desirable.

In another case two events, one random and one with an uncertain outcome, are jointly evaluated \cite{cheese}.
Both probabilistic and DS analyses yield likelihood estimates of the combined events equal to the probability
of the random event. This result is sometimes considered unsatisfactory, as the fusion process does not appear
to improve upon probability estimates of individual events \cite{cheese}. 
The result, however, is consistent with the frameworks of both analyses,
and presents no interpretation difficulties in a more general case, where the latter event is only partly uncertain.

This paper, by contrast, identifies a large class of bodies of evidence 
associated with non-zero conflict coefficient and yielding different DS and probabilistic evaluations that cannot be easily reconciled.  
The outcome sets are given by partitions and quasi-partitions of the set of evidence,
which correspond to the cases of zero and non-zero mass assignments to the universal set, respectively.
The finding contradicts a key result, an inequality that relates probabilistic and DS evaluations 
and thereby casts doubt on the legitimacy of probabilistic interpretation of the DS mass assignment 
when the DS rule of combination is used. 
This outcome supports the view expressed among others by Pearl that despite the formal similarity (and the precise algebraic relationship 
that was established between the two calculi in \cite{brodzik})
the DS and probabilistic approaches to evidence accumulation are separate theories with distinct objectives \cite{p2}, \cite{p3}.

\noindent
\section{Basic formulas}

\noindent
Denote by $\Omega$ a finite non-empty set of all possible outcomes of an event of interest, 
and by $2^{\Omega}$ the power set of $\Omega$. 
Define the set of observable outcomes, called the {\it set} (of subsets) {\it of evidence}, by 
\begin{eqnarray}
A=\{A_i \ | \ 0< i\leq |A| \} \subseteq 2^{\Omega}, \ \  A\neq \emptyset, \ \ 
\end{eqnarray}
where $|A|$ is the cardinality of $A$, and $\subseteq$ denotes "is a subset of". 

Given the set $A$ in (1), define a mapping 
\begin{equation}
m_A: 2^{\Omega}\mapsto [0,1], 
\end{equation}
such that
\begin{equation}
m_A(\emptyset)=0 
\end{equation}
and
\begin{equation}
\sum 
m_A(A_i)=1. 
\end{equation}
Set $m_{A_i}=m_A(A_i)$ and call it the {\it mass} of $A_i$. 
By an abuse of notation we will also write 
\begin{equation}
m_A=\{m_{A_i} \ | \ 0< i\leq |A| \}, 
\end{equation}
and refer to $m_A$ as the {\it mass assignment} of $A$.
Finally, we will call the set of pairs of the subsets $A_i$ and the corresponding masses $m_{A_i}$,
\begin{equation}
{\cal A}
=\{(A_i,m_{A_i}) \ | \ 0< i\leq |A|\},
\end{equation}
the {\it body of evidence} of $A$. 

The key difference between probability and mass is that probability is a measure and therefore 
it satisfies the additivity condition, that is, given a finite sequence $A_i$, $0<i\leq |A|$, of disjoint subsets of $A$,
\begin{equation}
P\left(\bigcup A_i\right)=\sum P(A_i).
\end{equation}
In general, mass does not satisfy condition (7).
Removing the additivity constraint can be convenient, as it permits inclusion of subjective judgments 
in the DS information fusion system, 
but it also has the undesirable consequence of making the interpretation of results of such fusion uncertain.
In particular, when considered together with the DS rule of combination, 
it is not always clear when mass can be made consistent with the standard probability evaluation.

Here we address this issue in a limited way by constraining mass to satisfy the additivity condition.
We identify mass with probability, combine bodies of evidence according to the DS rule, and test if mass
of the combined bodies agrees with the corresponding probabilities.
The additivity constraint imposed on mass allows us to focus on partitions and on bodies of evidence
with no contradictory mass assignments.
In the remainder of this section we explain the focus on partitions,
introduce the DS rule of combination, describe the auxiliary concepts of balance and plausibility, and
identify a key inequality linking probability and DS theories.

In general, $A$ may contain all non-trivial subsets of $2^{\Omega}$.
For example, when $A=\{a,b,c\}$, it is possible that $A=\{a,b,c,\{a,b\},\{a,c\},\{b,c\},\{a,b,c\}\}$.
Here, we restrict $A$ to be a {\it partition} of $\Omega$, i.e.,  
\begin{equation}
A_i\bigcap_{i\neq j} A_j=\emptyset 
\ \ {\rm and} \ \ \bigcup A_i=\Omega, 
\end{equation}
or a {\it quasi-partition} of $\Omega$, i.e.,
\begin{equation}
A_i\bigcap_{i\neq j\neq |A|}A_j=\emptyset, \ \ \bigcup_{i\neq |A|}A_i=\Omega \ \ {\rm and} \ \  A_{|A|}=\Omega.
\end{equation}
The latter case arises when the available information is uncertain, i.e., when $m_{\Omega}\neq 0$.
The reason for the restriction of sets of evidence to partitions 
is that it simplifies the analysis without removing generality:
provided condition (7) is satisfied,
bodies of evidence having overlapping sets can be replaced by bodies of evidence having no overlapping sets.
For example, the set $\{\{a,b\}, \{b,c\}\}$ can be replaced by the sets $\{\{a,b\},c\}$ and $\{a,\{b,c\}\}$.
Similarly, the set $\{\{a,b\}, \{a,b,c\},d\}$ can be replaced by the set $\{\{a,b\},c,d\}$. 

A key feature of the DS theory is the rule for combining bodies of evidence.
Let $\cal A$ and $\cal B$ be two distinct bodies of evidence.
Suppose a rule for combining the sets of evidence $A$ and $B$ and mapping the result to a decision set $C$, 
\begin{equation}
C=A\bigtriangledown B,
\end{equation}
is given by a partition of the set
\begin{equation}
\{A_i \cap B_j \ | \ 0< i\leq |A|, \ 0< j\leq |B| \}. 
\end{equation}
Assume an appropriate rule $\bigtriangledown$ is given.
The DS rule for combining the masses of $\cal A$ and $\cal B$ 
is then 
\begin{equation}
m_{C_k}
=
\frac{1}{1-\kappa}\sum_{A_i\cap B_j=
C_k}m_{A_i}m_{B_j}, \ \ 0<k\leq |C|,
\end{equation}
where
\begin{equation}
\kappa=\sum_{A_i\cap B_j=\emptyset}m_{A_i} m_{B_j}\neq 1\footnote{
In general, $0\leq \kappa\leq 1$.
$\kappa=1$ iff $\bigcup A_i \ \cap \ \bigcup B_j = \emptyset$, a satisfactory result,
since then $A$ and $B$ cannot be combined to form a decision set. 
For example, there might be bodies of evidence allowing one to evaluate outcomes 
"a tree is a poplar but not an oak" and "a tree is a cedar but not a pine",
but these cannot be combined to form a body of evidence allowing one to evaluate an outcome
"a tree is deciduous but not coniferous".} 
\end{equation}
is the {\it conflict coefficient} and 
\begin{equation}
{\cal C}
=\{(C_k,m_{C_k}) \ | \ 0< k\leq |C| \}
\end{equation}
is the DS composite body of evidence. 

Apart from mass, two other concepts are key in the DS theory: balance and plausibility.
{\it Balance} (or, {\it belief}) of a subset $A_i$ is the sum of the masses of all subsets $A_j$ of $A$, that are also subsets of
$A_i$, i.e.,
\begin{equation}
b_{A_i}=\sum_{A_j\subseteq A_i}m_{A_j}, \ \ 0<i\leq|A|. 
\end{equation}
{\it Plausibility} of a subset $A_i$ is the sum of the masses of all subsets $A_j$ of $A$, having non-empty intersection with $A_i$, i.e.,
\begin{equation}
p_{A_i}=\sum_{A_i\cap A_j\neq \emptyset}m_{A_j}, \ \ 0<i\leq |A|. 
\end{equation}
Like mass, balance and plausibility are mappings from the power set of $\Omega$ to the unit interval.
In particular, 
\begin{equation}
b_{\emptyset}=p_{\emptyset}=0
\end{equation}
and
\begin{equation}
b_{\Omega}=p_{\Omega}=1.
\end{equation}
Moreover, balance and plausibility are related by the formula
\begin{equation}
p_{A_i}=1-b_{\bar{A_i}}, \ \ 0<i\leq |A|, \ \ 
\bar{A_i}=\Omega - A_i.
\end{equation}
Due to Rota's generalization of the M\"obius inversion theorem \cite{rota}, mass can be uniquely recovered from balance by the formula
\begin{equation}
m_{A_j}=\sum_{A_i\subseteq A_j}(-1)^{|A_j-A_i|}b_{A_i}, \ \ 0<j\leq |A|.
\end{equation}
A similar formula exists for plausibility \cite{ds}; the two formulas ensure that no information is lost in
the process of performing (15) or (16).

A key result in DS theory describes the relationship among balance, plausibility and probability.
It follows from (15) and (16) that
\begin{equation}
b_{A_i}\leq p_{A_i}, \ \ 0< i\leq |A|.
\end{equation}
A stronger version of (21) that allows comparison of results of DS and probabilistic analyses has been proposed by Dempster
\cite{dempster} for the situation where mass assignment arises from a set-valued mapping from a probability space to $\Omega$,
\begin{equation}
b_{A_i}\leq P(A_i)\leq p_{A_i}, \ \ 0< i\leq |A|.
\end{equation}
Of particular importance to us are certain special cases. 
It follows from (8) and (9) that the condition (22) can be replaced by the condition
\begin{equation}
b_{A_i}= P(A_i)\leq p_{A_i}, \ \ 0< i\leq |A|,
\end{equation}
when $A$ is a quasi-partition of $\Omega$, and by the condition
\begin{equation}
b_{A_i}= P(A_i)= p_{A_i}, \ \ 0< i\leq |A|,
\end{equation}
when $A$ is a partition of $\Omega$. 
Since balance and plausibility bound the value of probability, 
they are often referred to as the {\it lower} and {\it upper probabilities}. 
A verification of validity of condition (22) and of its special cases, conditions (23) and (24), is the main goal of this paper. 

\noindent
\section{Combining bodies of evidence} 

\noindent
We analyze two cases of combining two bodies of evidence, both with a non-zero conflict coefficient, 
and both yielding inconsistent DS and probabilistic evaluations. In the first case the uncertainty mass of both bodies of evidence is zero. 
In the second case the uncertainty mass of one of the two bodies of evidence is non-zero. 
While the latter is a straightforward extension of the former, both cases are included for their pedagogical value.

\vspace{2mm}
 
\subsection{$m_{\Omega}=0$ and $\kappa\neq 0$}

\noindent
Consider the following two sets of evidence,
\begin{equation}
A
\doteq \{A_1, A_2 \}
=\{a,\{b,c\}\}
\end{equation}
and
\begin{equation}
B
\doteq \{B_1, B_2 \}
=\{\{a,b\},c\}, 
\end{equation}
having mass assignments
\begin{equation}
m_A\doteq \{m_{A_1}, m_{A_2}\}=\left\{\frac{1}{4}, \frac{3}{4}\right\}  
\end{equation}
and
\begin{equation}
m_B\doteq \{m_{B_1}, m_{B_2}\}=\left\{\frac{1}{2}, \frac{1}{2}\right\}. 
\end{equation}
Suppose the set combination rule is given by
\begin{eqnarray}
\nonumber
C
& = & A\bigtriangledown B \\
\nonumber
& \doteq & \{C_1=A_1\cap B_1, \ C_2=A_2\cap B_1, \  
 C_3=A_2\cap B_2 \} \\
& = & \{a, b, c\}.
\end{eqnarray}
We seek to obtain first, the mass of subsets of $C$, 
\begin{equation}
m_C\doteq \{m_{C_1}, m_{C_2}, m_{C_3}\}, 
\end{equation}
and second, the associated lower and upper probabilities. 

Since $A_1\cap B_2=\emptyset$, the conflict coefficient $\kappa=m_{A_1}m_{B_2}\neq 0$.
It follows from equation (12) that the mass of $C_1$, $C_2$ and $C_3$ is then
\begin{equation}
m_{C_1}
=\frac{m_{A_1}m_{B_1}}{1-m_{A_1}m_{B_2}} 
=\frac{\frac{1}{4}\frac{1}{2}}{1-\frac{1}{4}\frac{1}{2}}
=\frac{1}{7}, 
\end{equation} 
\begin{equation}
m_{C_2}
=\frac{m_{A_2}m_{B_1}}{1-m_{A_1}m_{B_2}} 
=\frac{\frac{3}{4}\frac{1}{2}}{1-\frac{1}{4}\frac{1}{2}}
=\frac{3}{7} 
\end{equation}
and
\begin{equation}
m_{C_3}
=\frac{m_{A_2}m_{B_2}}{1-m_{A_1}m_{B_2}} 
=\frac{\frac{3}{4}\frac{1}{2}}{1-\frac{1}{4}\frac{1}{2}}
=\frac{3}{7}. 
\end{equation} \\
Since $C$ is a partition, then $m_{C_i}=b_{C_i}=p_{C_i}, \ i=1,2,3$,
and we are done. 

Suppose the mass assignments (27) and (28) coincide with probabilities. We will treat these two mass assignments as partial
information about a fixed probability distribution that we seek to derive. It follows then, that
\begin{equation}
P(C_1)=
P(A_1)=\frac{1}{4},
\end{equation} 
\begin{equation}
P(C_2)=
P(A_2)-P(B_2)=\frac{3}{4}-\frac{1}{2}=\frac{1}{4}
\end{equation}
and
\begin{equation}
P(C_3)=
P(B_2)=\frac{1}{2}.
\end{equation} \\
Comparing rhs of equations (31)-(33) and (34)-(36), we have 
\begin{equation}
m_{C_i}\neq P(C_i), \ \ i=1,2,3,  
\end{equation} 
and therefore the condition (24) is not satisfied.

To verify if the result in (37) is an anomaly, consider a general case, given by the mass assignment 
\begin{equation}
m_A=\{x, 1-x\} 
\end{equation}
and
\begin{equation}
m_B=\{y, 1-y\}, 
\end{equation}
$0\leq x,y\leq 1$.
Then from equation (12)
\begin{equation}
m_{C_1}
=\frac{m_{A_1}m_{B_1}}{1-m_{A_1}m_{B_2}} 
=\frac{xy}{1-x(1-y)}, 
\end{equation} 
\begin{equation}
m_{C_2}
=\frac{m_{A_2}m_{B_1}}{1-m_{A_1}m_{B_2}} 
=\frac{(1-x)y}{1-x(1-y)}
\end{equation}
and
\begin{equation}
m_{C_3}
=\frac{m_{A_2}m_{B_2}}{1-m_{A_1}m_{B_2}} 
=\frac{(1-x)(1-y)}{1-x(1-y)}. 
\end{equation} 
As before, suppose the mass assignment (38)-(39) coincides with probabilities.
It then follows from (25)-(26) and (38)-(39) that 
\begin{equation}
P(C_1)=
P(A_1)=x,
\end{equation}
\begin{equation}
P(C_2)=
P(B_1)-P(A_1)=y-x
\end{equation}
and
\begin{equation}
P(C_3)=
P(B_2)=1-y.
\end{equation}
Comparing rhs of equations (40)-(42) and (43)-(45), it follows that mass and probabilities are equal and therefore the condition (24) is satisfied 
if and only if
$x=0$ and $y$ is arbitrary, or $y=1$ and $x$ is arbitrary. 
This condition is equivalent to the condition $\kappa=0$. 

\vspace{2mm}

\noindent
\subsection{$m_{\Omega} \ne 0$ and $\kappa \ne 0$}  

\noindent
Consider the following two sets of evidence,
\begin{equation}
A
\doteq \{A_1, A_2, A_3 \}
=\{a,\{b,c\}, \{a,b,c\}\}
\end{equation}
and
\begin{equation}
B
\doteq \{B_1, B_2 \}
=\{\{a,b\},c\}, 
\end{equation}
having mass assignments
\begin{equation}
m_A\doteq \{m_{A_1}, m_{A_2}, m_{A_3}\}=\left\{x, \bar{x}, 1-x-\bar{x}\right\}  
\end{equation}
and
\begin{equation}
m_B\doteq \{m_{B_1}, m_{B_2}\}=\left\{y, 1-y\right\}, 
\end{equation}
where $0\leq x+\bar{x}, y\leq 1$.
Suppose the set combination rule is given by
\begin{eqnarray}
\nonumber
C & = & A\bigtriangledown B \\
\nonumber
& \doteq & \{C_1=A_1\cap B_1, \ C_2=A_2\cap B_1, \  \\
\nonumber
& & C_3=A_2\cap B_2 \ \cup \ A_3\cap B_2, \ C_4=A_3\cap B_1 \} \\
& = & \{a, b, c, \{a,b\}\}.
\end{eqnarray}
We seek to obtain
\begin{equation}
m_C\doteq \{m_{C_1}, m_{C_2}, m_{C_3}, m_{C_4}\}, 
\end{equation}
and the associated values of balance and plausibility.
Since $A_1\cap B_2=\emptyset$, the conflict coefficient $\kappa=m_{A_1}m_{B_2}=x(1-y)\neq 0$, except in the trivial case.
It follows from equation (12) that the masses of $C_1$, $C_2$, $C_3$ and $C_4$ are then
\begin{equation}
m_{C_1}
=\frac{m_{A_1}m_{B_1}}{1-m_{A_1}m_{B_2}} 
=\frac{xy}{1-x(1-y)},
\end{equation} 
\begin{equation}
m_{C_2}
=\frac{m_{A_2}m_{B_1}}{1-m_{A_1}m_{B_2}} 
=\frac{\bar{x}y}{1-x(1-y)},
\end{equation}
\begin{equation}
m_{C_3}
=\frac{m_{A_2}m_{B_2}+m_{A_3}m_{B_2}}{1-m_{A_1}m_{B_2}} 
=\frac{(1-x)(1-y)}{1-x(1-y)},
\end{equation}
and
\begin{equation}
m_{C_4}
=\frac{m_{A_3}m_{B_1}}{1-m_{A_1}m_{B_2}} 
=\frac{(1-x-\bar{x})y}{1-x(1-y)},
\end{equation} \\
where $x(1-y)\neq 1$, 
and, from equations (15)-(16) that the corresponding balances and plausibilities are, respectively
\begin{equation}
b_C\doteq \{b_{C_1}, b_{C_2}, b_{C_3}\}=\left\{m_{C_1}, m_{C_2}, m_{C_3}\right\}  
\end{equation}
and
\begin{equation}
p_C\doteq \{p_{C_1}, p_{C_2}, p_{C_3}\}=\left\{m_{C_1}+m_{C_4}, m_{C_2}+m_{C_4}, m_{C_3}\right\}. 
\end{equation}

The objective, as before, is to evaluate the consistency of results generated by identification of DS masses with probabilities.
Equating the probability of $c$ with the mass of $B_2$, we have
\begin{equation}
P(c)=1-y.
\end{equation}
Furthermore, since 
by (56) and (57), 
\begin{equation}
m_{C_3}=b_{C_3}=p_{C_3}, 
\end{equation}
then, by (22) and (50), 
\begin{equation}
m_{C_3}=P(C_3)=P(c).
\end{equation}
Combining the last two results yields 
\begin{equation}
1-y=\frac{(1-x)(1-y)}{1-x(1-y)}.
\end{equation}
Equation (61) is satisfied if and only if $x=0$ or $y=0$ or $y=1$.
The first and the last case imply $\kappa=0$. The second case implies $\kappa=x$.
However, since $P(c)=1$ then $P(a)=0$ and therefore, as before, $x=0$. 

Similarly inconsistent evaluations are obtained for the singletons ${a}$ and ${b}$. 
The evaluations of $P(b)$ for $x=1/4$, $\bar{x}=1/2$ and $y=1/2$ are particularly revealing. 
Substituting $x$ and $y$ in (53) and (55) and proceeding as before leads to the DS evaluation
\begin{equation}
\frac{2}{7}\leq P(b) \leq \frac{3}{7}
\end{equation}
and the probability evaluation 
\begin{equation}
0 \leq P(b) \leq \frac{1}{4}.
\end{equation}
Note that the two evaluations are not merely different - they do not overlap! 
This anomaly cannot be relieved by renormalization of balance and plausibility 
suggested in \cite{bp}; in fact the problem then becomes even more severe.

It follows from the preceding argument that given two bodies of evidence 
equipped with an arbitrary mass assignment and an arbitrary set combination rule,
but satisfying the non-zero conflict coefficient condition, use of the DS rule of combination can yield a 
mass assignment for the combined body of evidence that is inconsistent with probabilities, thereby violating the inequality (22).

\end{document}